\documentclass[graybox]{svmult}

\usepackage{mathptmx}       
\usepackage{helvet}         
\usepackage{courier}        
\usepackage{type1cm}        
\usepackage{makeidx}         
\usepackage{graphicx}        
\usepackage{multicol}        
\usepackage[bottom]{footmisc}
\usepackage{amsfonts}
\usepackage{bm}
\usepackage{txfonts}
\usepackage{color}

\makeindex             

\begin{document}

\title*{Algebraic dual polynomials for the equivalence of curl-curl problems}
\author{Marc Gerritsma, Varun Jain, Yi Zhang  and Artur Palha}
\institute{Marc Gerritsma, Varun Jain and Yi Zhang \at Faculty of Aerospace Engineering, TU Delft, Kluyverweg 1, 2629 HS, Delft, \email{\{M.I.Gerritsma,V.Jain,Y.Zhang-14\}@tudelft.nl}
\and Artur Palha \at Department of Mechanical Engineering, Technical University Eindhoven, \email{A.PalhaDaSilvaClerigo@tudelft.nl}}
%
%
\maketitle

\abstract*{In this paper we will consider two curl-curl equation in two dimensions. One curl-curl problem for a scalar quantity $F$ and one problem for a vector field $\bm{E}$. For Dirichlet boundary conditions $\bm{n} \times \bm{E} = \hat{E}_{\dashv}$ on $\bm{E}$ and Neumann boundary conditions $\bm{n} \times \mbox{\bf{curl}}\,F=\hat{E}_{\dashv}$, we expect the solutions to satisfy $\bm{E}=\mbox{\bf{curl}}\,F$. When we use algebraic dual polynomial representations, these identities continue to hold at the discrete level. Equivalence will be proved and illustrated with a computational example.}

\abstract{In this paper we will consider two curl-curl equation in two dimensions. One curl-curl problem for a scalar quantity $F$ and one problem for a vector field $\bm{E}$. For Dirichlet boundary conditions $\bm{n} \times \bm{E} = \hat{E}_{\dashv}$ on $\bm{E}$ and Neumann boundary conditions $\bm{n} \times \mbox{\bf{curl}}\,F=\hat{E}_{\dashv}$, we expect the solutions to satisfy $\bm{E}=\mbox{\bf{curl}}\,F$. When we use algebraic dual polynomial representations, these identities continue to hold at the discrete level. Equivalence will be proved and illustrated with a computational example.}

\section{Introduction}\label{sec:intro}
Numerical methods lead invariably to approximations, but a judicious choice of finite dimensional function spaces allows one to preserve, at the discrete level, identities that hold at the continuous level. In this paper we will focus on the finite dimensional representation of the curl operator; or, to put it more correctly, the {\em curl operators}. This is particularly clear in the two-dimensional setting where one curl operator maps scalar fields to vector fields and the other curl operator is its adjoint and therefore maps vector fields to scalar fields.

Faraday's and Amp\`{e}re's law demonstrate the importance of the curl operator in electromagnetism. In fluid mechanics the curl operator appears in the relation between the stream function and the mass fluxes, and the definition of vorticity, \cite{MEEVC}.

In $\mathbb{R}^2$ we define the curl operators $\mbox{\bf{curl}}\,F = \left ( \partial F/\partial y,-\partial F/\partial x \right )^{\top}$ for a scalar field $F(x,y)$ and $\mbox{curl}\,\bm{E} = \partial E_y/\partial x - \partial E_x/\partial y$ for a vector field $\bm{E}=\left ( E_x,E_y \right )$. We define the functions spaces
\begin{equation}
H(\mbox{curl};\Omega) = \left \{ \bm{E} \in \left [ L^2(\Omega) \right ]^2 \,\big | \, \mbox{curl}\,\bm{E} \in L^2(\Omega) \, \right \} \;\; \mbox{and} \;,
\end{equation}
and
\begin{equation}
H(\mbox{\bf{curl}};\Omega) = \left \{ F \in  L^2(\Omega) \,\big | \, \mbox{\bf{curl}}\,F \in \left [ L^2(\Omega) \right ]^2 \, \right \} \;.
\end{equation}

\section{The equivalent curl-curl dual problems}\label{sec:curl_curl_problem}
In \cite{Carstensen} the equivalence between several Dirichlet and Neumann problems is introduced. The main question we want to address in this paper is whether we can preserve these equivalences at the discrete level. In \cite{Jain} this equivalence was already established at the discrete level for the scalar grad-div problem. In this paper we want to focus on the curl-curl equivalence problem: Given $\hat{E}_{\dashv} \in H^{-\frac{1}{2}}(\mbox{div};\partial \Omega)$, \cite{BuffaCiarlet}, find the solution $\bm{E} \in H(\mbox{curl};\Omega)$ of the Dirichlet problem satisfying
\begin{equation}
\left \{ \begin{array}{ll}
n \times \bm{E} = \hat{E}_{\dashv} \quad \quad &\mbox{ on } \partial \Omega \\[1ex]
\mbox{\bf{curl}} \left ( \mbox{curl}\, \bm{E} \right ) + \bm{E} = \bm{0} \quad \quad &\mbox{in } \Omega
\end{array} \right . \; ,
\label{eq:curlE}
\end{equation}
and the associated Neumann problem given by: For $\hat{E}_{\dashv} \in H^{-\frac{1}{2}}(\mbox{div};\partial \Omega)$ find the solution $F \in H(\mbox{\bf{curl}};\Omega)$ such that
\begin{equation}
\left \{ \begin{array}{ll}
n \times \mbox{\bf{curl}}\,F = \hat{E}_{\dashv} \quad \quad &\mbox{ on } \partial \Omega \\[1ex]
\mbox{curl} \left ( \mbox{\bf{curl}} \,F \right ) + F = 0 \quad \quad &\mbox{in } \Omega
\end{array} \right . \; .
\label{eq:curlF}
\end{equation}
At the continuous level we know that these two problems are equivalent in the sense that $\bm{E}$ solves (\ref{eq:curlE}) and $F$ solves (\ref{eq:curlF}) if and only if $\bm{E}=\mbox{\bf{curl}}\,F$. In addition, we have that $\|\bm{E}\|_{H(\mbox{curl};\Omega)}=\|F\|_{H(\mbox{\bf{curl}};\Omega)} = \| \hat{E}_{\dashv} \|_{H^{-\frac{1}{2}}(\mbox{div};\partial \Omega)}$.

In this paper we want to present a spectral element formulation for both problems (\ref{eq:curlE}) and (\ref{eq:curlF}) such that the equivalence, $\bm{E}^h=\mbox{\bf{curl}}\,F^h$, continues to hold for the discrete solutions. Moreover, we want to show that we have $\|\bm{E}^h\|_{H(\mbox{curl};\Omega)}=\|F^h\|_{H(\mbox{\bf{curl}};\Omega)}$.

\section{Primal spectral element formulation}
\label{sec:primal_formulation}
Consider the partitioning $-1=\xi_0<\xi_1 <,\ldots,<\xi_{N-1}<\xi_N=1$ of the interval $I=[-1,1]$, where $\xi_i$ are the roots of the $(1-\xi^2)L_N^{\prime}(\xi)$, with $L_N^{\prime}(\xi)$ the derivative of the Legendre polynomial of degree $N$. With these nodes $\xi_i$ we associate the Lagrange polynomials, $h_i(\xi)$, of degree $N$ which satisfy $h_i(\xi_j)=\delta_{ij}$. Any polynomial $p$ of degree $N$ defined on $I$ can be written as
\begin{equation}
p(\xi) = \sum_{i=0}^N p_i h_i(\xi) \;.
\label{eq:1D_nodal_expansion}
\end{equation}
Since the Lagrange polynomials are linearly independent, $p(\xi)=0$ iff all $p_i=0$.

The derivative of (\ref{eq:1D_nodal_expansion}) is given by
\begin{equation}
p^{\prime}(\xi) = \sum_{i=0}^N p_i \frac{dh_i}{d\xi} (\xi) = \sum_{i=1}^N \left ( p_i - p_{i-1} \right ) e_i(\xi) \;,
\label{eq:1D_deriivative}
\end{equation}
where, \cite{Gerritsma},
\[ e_i(\xi) = - \sum_{k=1}^{i-1} \frac{dh_k}{d\xi}(\xi) \;. \]
Note that the $dh_i/d\xi$, $i=0,\ldots,N$, do not form a basis, while the functions $e_i(\xi)$, $i=1,\ldots,N$ do form a basis. Therefore, $p^{\prime}(\xi) =0$ iff $p_i=p_{i-1}$ for $i=1,\ldots,N$, which in turn means that $p(\xi)=const$, as required. Note that differentiation of the Lagrange expansion (\ref{eq:1D_nodal_expansion}) amounts to a linear combination of the expansion coefficients, $(p_i-p_{i-1})$ in (\ref{eq:1D_deriivative}) and a representation in a different basis, $e_i(\xi)$ in (\ref{eq:1D_deriivative}). An important property of the edge polynomials, $e_i(\xi)$, is
\begin{equation}
\int_{\xi_{j-1}}^{\xi_j} e_i(\xi)\,\mathrm{d}\xi = \left \{ \begin{array}{ll}
1 \quad \quad & \mbox{if } i=j \\[1ex]
0 \quad \quad & \mbox{if } i\neq j
\end{array} \right . \;.
\label{eq:integral_property_edge}
\end{equation}

\begin{center}
\begin{figure}[h!]
\includegraphics[scale=.55]{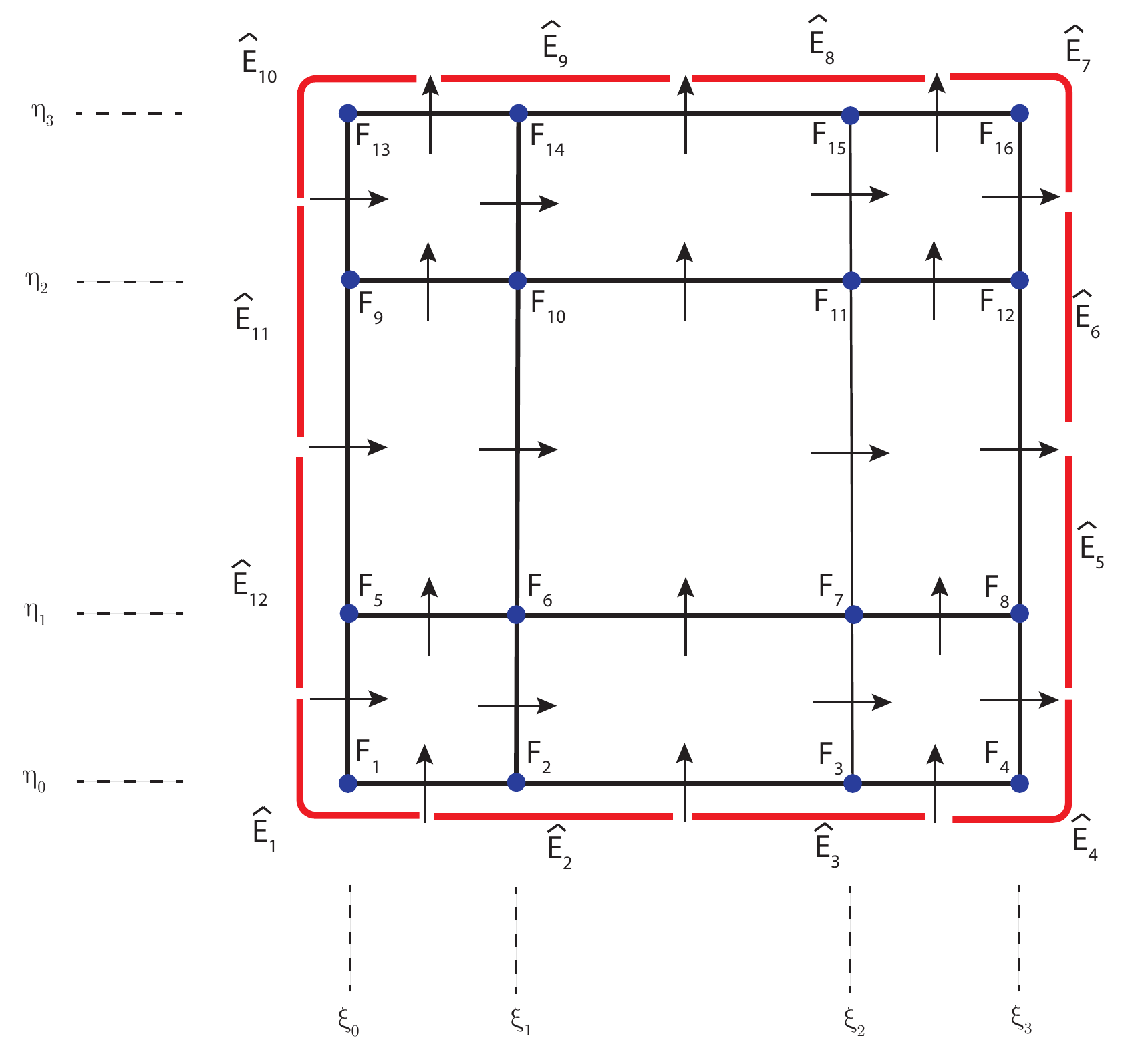}
%
%
\caption{Layout of one spectral element: The degrees of freedom for $F$ indicated by the blue points. The curl of $F$ is represented on the mesh line segments by the arrows crossing the line segments in the mesh. The boundary condition $\bm{n}\times \mbox{\bf{curl}}\,F=\hat{E}_{\dashv}$ is represented by the red outer line segments.}
\label{fig:F_grid}       
\end{figure}
\end{center}

In the two dimensional case we consider $(\xi,\eta) \in \hat{\Omega} = [-1,1]^2\subset \mathbb{R}^2$ and the partitioning $\xi_i$ in the $\xi$-direction as given in the one dimensional case and we choose the same partitioning in the $\eta$-direction, see Figure~\ref{fig:F_grid}. Here $\Omega$ is a contractible domain with Lipschitz continuous boundary. For the representation of $F$ we use the tensor product of the nodal representation
\begin{equation}
F^h(\xi,\eta) = \sum_{i=0}^N \sum_{j=0}^N F_{i,j}h_i(\xi)h_j(\eta) = \Psi^{(0)}(\xi,\eta) \mathsf{F} \;.
\label{eq:F_expansion}
\end{equation}
Here $\Psi^{(0)}(\xi,\eta)$ is a row vector with the basis functions and $\mathsf{F}$ is a column vector with the nodal degrees of freedom
\begin{equation}
\Psi^{(0)}(\xi,\eta) = \left [ \begin{array}{ccc}
h_0(\xi)h_0(\eta) \quad  & \ldots \quad  & h_N(\xi)h_N(\eta)
\end{array} \right ] \;,\;\;\; \mathsf{F} = \left [ \begin{array}{c}
F_{0,0}\\
F_{1,0}\\
\vdots \\
F_{N-1,N}\\
F_{N,N}
\end{array} \right ] \;,
\label{eq:primal_0_basis}
\end{equation}
where $F_{i,j}=F^h(\xi_i,\eta_j)$.
The curl of $F^h$ is then given by
\begin{equation}
\mbox{\bf{curl}}\,F^h = \left ( \begin{array}{c}
\sum_{i=0}^N \sum_{j=1}^N (F_{i,j}-F_{i,j-1})h_i(\xi)e_j(\eta) \\[1ex]
\sum_{i=1}^N \sum_{j=0}^N (-F_{i,j}+F_{i-1,j})e_i(\xi)h_j(\eta)
\end{array} \right ) = \Psi^{(1)}(\xi,\eta) \mathbb{E}^{1,0} \mathsf{F} \;.
\label{eq:curlF_expansion}
\end{equation}
In (\ref{eq:curlF_expansion}) we used (\ref{eq:1D_deriivative}) to represent the curl in a different basis. The basis $\Psi^{(1)}(\xi,\eta)$ and the incidence matrix $\mathbb{E}^{1,0}$ are given by
\begin{equation}
\Psi^{(1)}(\xi,\eta) = \left [ \begin{array}{cccccc}
h_0(\xi)e_1(\eta) & \ldots & h_N(\xi)e_N(\eta) & 0 & \ldots & 0 \\[1ex]
0 & \ldots & 0 & e_1(\xi)h_0(\eta) & \ldots & e_N(\xi)h_N(\eta)
\end{array} \right ] \;,
\label{eq:primal_1_basis}
\end{equation}
\[
\mathbb{E}^{1,0} = \left ( \begin{array}{cccccccccccccccc}
-1 & 0 & 0 & 0 & 1 & 0 & 0 & 0 & 0 & 0 & 0 & 0 & 0 & 0 & 0 & 0 \\
0 & -1 & 0 & 0 & 0 & 1 & 0 & 0 & 0 & 0 & 0 & 0 & 0 & 0 & 0 & 0 \\
0 & 0 & -1 & 0 & 0 & 0 & 1 & 0 & 0 & 0 & 0 & 0 & 0 & 0 & 0 & 0 \\
0 & 0 & 0 & -1 & 0 & 0 & 0 & 1 & 0 & 0 & 0 & 0 & 0 & 0 & 0 & 0 \\
0 & 0 & 0 & 0 & -1 & 0 & 0 & 0 & 1 & 0 & 0 & 0 & 0 & 0 & 0 & 0 \\
0 & 0 & 0 & 0 & 0 & -1 & 0 & 0 & 0 & 1 & 0 & 0 & 0 & 0 & 0 & 0 \\
0 & 0 & 0 & 0 & 0 & 0 & -1 & 0 & 0 & 0 & 1 & 0 & 0 & 0 & 0 & 0 \\
0 & 0 & 0 & 0 & 0 & 0 & 0 & -1 & 0 & 0 & 0 & 1 & 0 & 0 & 0 & 0 \\
0 & 0 & 0 & 0 & 0 & 0 & 0 & 0 & -1 & 0 & 0 & 0 & 1 & 0 & 0 & 0 \\
0 & 0 & 0 & 0 & 0 & 0 & 0 & 0 & 0 & -1 & 0 & 0 & 0 & 1 & 0 & 0 \\
0 & 0 & 0 & 0 & 0 & 0 & 0 & 0 & 0 & 0 & -1 & 0 & 0 & 0 & 1 & 0 \\
0 & 0 & 0 & 0 & 0 & 0 & 0 & 0 & 0 & 0 & 0 & -1 & 0 & 0 & 0 & 1 \\
1 & -1 & 0 & 0 & 0 & 0 & 0 & 0 & 0 & 0 & 0 & 0 & 0 & 0 & 0 & 0 \\
0 & 1 & -1 & 0 & 0 & 0 & 0 & 0 & 0 & 0 & 0 & 0 & 0 & 0 & 0 & 0 \\
0 & 0 & 1 & -1 & 0 & 0 & 0 & 0 & 0 & 0 & 0 & 0 & 0 & 0 & 0 & 0 \\
0 & 0 & 0 & 0 & 1 & -1 & 0 & 0 & 0 & 0 & 0 & 0 & 0 & 0 & 0 & 0 \\
0 & 0 & 0 & 0 & 0 & 1 & -1 & 0 & 0 & 0 & 0 & 0 & 0 & 0 & 0 & 0 \\
0 & 0 & 0 & 0 & 0 & 0 & 1 & -1 & 0 & 0 & 0 & 0 & 0 & 0 & 0 & 0 \\
0 & 0 & 0 & 0 & 0 & 0 & 0 & 0 & 1 & -1 & 0 & 0 & 0 & 0 & 0 & 0 \\
0 & 0 & 0 & 0 & 0 & 0 & 0 & 0 & 0 & 1 & -1 & 0 & 0 & 0 & 0 & 0 \\
0 & 0 & 0 & 0 & 0 & 0 & 0 & 0 & 0 & 0 & 1 & -1 & 0 & 0 & 0 & 0 \\
0 & 0 & 0 & 0 & 0 & 0 & 0 & 0 & 0 & 0 & 0 & 0 & 1 & -1 & 0 & 0 \\
0 & 0 & 0 & 0 & 0 & 0 & 0 & 0 & 0 & 0 & 0 & 0 & 0 & 1 & -1 & 0 \\
0 & 0 & 0 & 0 & 0 & 0 & 0 & 0 & 0 & 0 & 0 & 0 & 0 & 0 & 1 & -1
\end{array} \right ) \;,
\]
where this incidence matrix corresponds to the layout depicted in Figure~\ref{fig:F_grid}, i.e. $N=3$. Note that this incidence matrix only contains the entries $-1$, $0$ and $1$ and that the matrix is extremely sparse. The important thing to note is that this incidence matrix remains unchanged if we map the standard element $\hat{\Omega}$ to an arbitrary curved element. The basis functions $\Psi^{(1)}(\xi,\eta)$ do change, but the incidence matrix remains invariant. This is another reason to decompose a derivative into a part that acts on the degrees of freedom and new basis functions, as was done in (\ref{eq:1D_deriivative}).

The mass matrix $\mathbb{M}^{(0)}$ associated with the basis functions (\ref{eq:primal_0_basis}) is given by
\begin{equation}
\mathbb{M}^{(0)} = \int_{\hat{\Omega}} {\Psi^{(0)}}^{\top} \Psi^{(0)}\,\mathrm{d} \Omega \;.
\label{eq:M0}
\end{equation}

Likewise, the mass matrix $\mathbb{M}^{(1)}$ associated with the basis functions (\ref{eq:primal_1_basis}) is given by
\begin{equation}
\mathbb{M}^{(1)} = \int_{\hat{\Omega}} {\Psi^{(1)}}^{\top} \Psi^{(1)}\,\mathrm{d} \Omega \;.
\label{eq:M1}
\end{equation}

\section{Dual spectral element formulation}\label{sec:dual_formulation}
\subsection{Duality in the interior of the domain}
In the previous section we expanded the discrete solution in terms of basis functions $h_i(\xi)h_j(\eta)$ for $F^h$ and $h_i(\xi)e_j(\eta) \otimes e_i(\xi)h_j(\eta)$ for the curl of $F^h$, respectively. With every linear vector space, $\mathcal{V}$, we can associate the space of linear functionals acting on that space $\mathcal{L}(\mathcal{V},\mathbb{R})=\mathcal{V}^{\prime}$, called the {\em algebraic dual space}. Let $\alpha \in \mathcal{V}^{\prime}$ and $u\in \mathcal{V}$, then $\alpha(u)\in \mathbb{R}$. Because we work in a Hilbert space, the Riesz representation theorem tells us that for every $\alpha \in \mathcal{V}^{\prime}$ there exists a unique $v_{\alpha} \in \mathcal{V}$ such that
\begin{equation}
\alpha(u) = (v_{\alpha},u)\;,\;\;\; \forall u \in \mathcal{V} \;,
\label{eq:Riesz}
\end{equation}
where $(\cdot,\cdot)$ denotes the inner product in $\mathcal{V}$, \cite{Kreyszig,OdenDemkowicz}. We first apply these ideas to the degrees of freedom (the expansion coefficients) which also form a linear vector space. Let $F^h$ and $G^h$ be expanded as in (\ref{eq:F_expansion})
\[ F^h(\xi,\eta) = \Psi^{(0)}(\xi,\eta) \mathsf{F}\quad \mbox{ and } \quad  G^h(\xi,\eta) = \Psi^{(0)}(\xi,\eta) \mathsf{G} \;.\]
Then we define {\em the dual degrees of freedom} $\widetilde{\mathsf{F}}$ analogous to (\ref{eq:Riesz}) by, \cite{Jain,Yi}
\begin{equation}
\widetilde{\mathsf{F}}^{\top} \mathsf{G}:= \mathsf{F}^{\top} \mathbb{M}^{(0)} \mathsf{G} \;, \;\;\; \forall \mathsf{G} \in \mathbb{R}^{(N+1)^2} \;.
\end{equation}
Therefore, the dual degrees of freedom are related to the primal degrees of freedom by $\widetilde{\mathsf{F}} = \mathbb{M}^{(0)} \mathsf{F}$. The canonical dual basis functions are then given by
\begin{equation}
\widetilde{\Psi}^{(2)} := \Psi^{(0)} {\mathbb{M}^{(0)}}^{-1}\;,
\label{eq:dual_volume_basis}
\end{equation}
such that
\begin{equation}
\int_{\Omega} {\widetilde{\Psi}^{(2)^{\top}}} \Psi^{(0)} \,\mathrm{d}\Omega = \mathbb{I} \;,
\label{eq:canonica_dual}
\end{equation}
where $\mathbb{I}$ is the identity matrix on $\mathbb{R}^{(N+1)^2}$. The relation (\ref{eq:canonica_dual}) is analogous to the canonical basis $e^*_i \in \mathcal{V}^{\prime}$ with the property $e^*_i(e_j) = \delta_{ij}$, when $e_i$ form a basis for $\mathcal{V}$. If the basis functions change under a transformation, then the dual basis functions also change and the (\ref{eq:canonica_dual}) continues to hold.

Let the vector field $\bm{E}^h$ be expanded as in (\ref{eq:curlF_expansion})
\[ \bm{E}^h(\xi,\eta) = \Psi^{(1)}(\xi,\eta) \mathsf{E} \;.\]
The corresponding dual degrees of freedom are then given by $\widetilde{\mathsf{E}} = \mathbb{M}^{(1)}\mathsf{E}$ and the associated dual basis is related to the primal basis by $\widetilde{\Psi}^{(1)}(\xi,\eta) = \Psi^{(1)}(\xi,\eta) {\mathbb{M}^{(1)}}^{-1}$.

\subsection{Duality in the boundary}
\label{sec:duality_boundary}
The construction of a primal and a dual representation in the interior of the domain $\hat{\Omega}$ can also be applied along the boundary of the domain $\partial \hat{\Omega}$. We can restrict $F^h$ to the boundary of the domain using (\ref{eq:F_expansion}), which gives
\begin{equation}
 \left . F^h \right |_{\partial \hat{\Omega}} = \left \{ \begin{array}{l}
 \left . F^h \right |_{\partial \hat{\Omega}_E} = F_{N,j}h_j(\eta) \\[1ex]
 \left . F^h \right |_{\partial \hat{\Omega}_N} = F_{i,N}h_i(\xi) \\[1ex]
 \left . F^h \right |_{\partial \hat{\Omega}_W} = F_{0,j}h_j(\eta) \\[1ex]
 \left . F^h \right |_{\partial \hat{\Omega}_S} = F_{i,0}h_j(\xi)
 \end{array} \right . \;.
\end{equation}
This boundary expansion is essentially a one-dimensional expansion, (\ref{eq:1D_nodal_expansion}), in terms of the four 1D elements which make up the boundary of a single spectral element. From this expansion we can compute the associated mass matrix, which for this boundary integral we will denote by $\mathbb{B}^{(0)}$. With the nodal degrees of freedom on the boundary, $\mathsf{F}_b$, we can now define the dual degrees of freedom, $\widetilde{\mathsf{E}}_b$ by setting $\widetilde{\mathsf{E}}_b = \mathbb{B}^{(0)} \mathsf{F}_b$.

We introduce the matrix $\mathbb{T}$, given by
\begin{equation}
\mathbb{T} = \left ( \begin{array}{cccccccccccccccc}
1 & 0 & 0 & 0 & 0 & 0 & 0 & 0 & 0 & 0 & 0 & 0 & 0 & 0 & 0 & 0 \\
0 & 1 & 0 & 0 & 0 & 0 & 0 & 0 & 0 & 0 & 0 & 0 & 0 & 0 & 0 & 0 \\
0 & 0 & 1 & 0 & 0 & 0 & 0 & 0 & 0 & 0 & 0 & 0 & 0 & 0 & 0 & 0 \\
0 & 0 & 0 & 1 & 0 & 0 & 0 & 0 & 0 & 0 & 0 & 0 & 0 & 0 & 0 & 0 \\
0 & 0 & 0 & 0 & 0 & 0 & 0 & 1 & 0 & 0 & 0 & 0 & 0 & 0 & 0 & 0 \\
0 & 0 & 0 & 0 & 0 & 0 & 0 & 0 & 0 & 0 & 0 & 1 & 0 & 0 & 0 & 0 \\
0 & 0 & 0 & 0 & 0 & 0 & 0 & 0 & 0 & 0 & 0 & 0 & 0 & 0 & 0 & 1 \\
0 & 0 & 0 & 0 & 0 & 0 & 0 & 0 & 0 & 0 & 0 & 0 & 0 & 0 & 1 & 0 \\
0 & 0 & 0 & 0 & 0 & 0 & 0 & 0 & 0 & 0 & 0 & 0 & 0 & 1 & 0 & 0 \\
0 & 0 & 0 & 0 & 0 & 0 & 0 & 0 & 0 & 0 & 0 & 0 & 1 & 0 & 0 & 0 \\
0 & 0 & 0 & 0 & 0 & 0 & 0 & 0 & 1 & 0 & 0 & 0 & 0 & 0 & 0 & 0 \\
0 & 0 & 0 & 0 & 1 & 0 & 0 & 0 & 0 & 0 & 0 & 0 & 0 & 0 & 0 & 0
\end{array} \right ) \;.
\end{equation}
The matrix $\mathbb{T}$ restricts the field $F^h$ to the boundary of the domain.
The curl of the representation $\widetilde{\bm{E}}^h$ is defined in the weak sense as
\begin{eqnarray}
\int_{\hat{\Omega}} F^h \cdot \mbox{curl}\,\widetilde{\bm{E}}^h\,\mathrm{d} \Omega & = & \int_{\hat{\Omega}} \mbox{\bf{curl}} F^h \cdot \widetilde{\bm{E}}^h \,\mathrm{d} \Omega + \int_{\partial \hat{\Omega}}  F^h \cdot \bm{n} \times \widetilde{\bm{E}}^h \,\mathrm{d} \Gamma \nonumber \\[1ex]
 & = & \mathsf{F}^{\top} {\mathbb{E}^{1,0}}^{\top} \widetilde{\mathsf{E}} + {\mathsf{F}}^{\top} \mathbb{T}^{\top} \widetilde{\mathsf{E}}_b  \;,\;\;\; \forall \mathsf{F} \;,\label{eq:dirivative_dual_E}
\end{eqnarray}
where $\widetilde{\mathsf{E}}_b$ are the degrees of freedom along the boundary indicated by the red line segments in Fig.~\ref{fig:F_grid}. Note also that minus signs in $
\bm{n}\times \bm{E}$ cancel with the minus signs originating from the counter-clockwise evaluated boundary integral in (\ref{eq:dirivative_dual_E}).
Also, (\ref{eq:dirivative_dual_E}) shows that the degrees of freedom for curl of $\widetilde{\bm{E}}^h$ are given by ${\mathbb{E}^{1,0}}^{\top} \widetilde{\mathsf{E}} + \mathbb{T}^{\top} \widetilde{\mathsf{E}}_b$.

\section{Discrete formulation of the curl-curl problem}\label{sec:Discrete_formulations_curl_curl}
\subsection{The Neumann problem}
The variational form of the Neumann problem (\ref{eq:curlF}) is given by: Find $F^h \in H(\mbox{\bf{curl}};\Omega)$ such that
\begin{equation} \label{eq:var_neum}
(\mbox{\bf{curl}}\ G^h, \mbox{\bf{curl}}\ F^h) + (G^h, F^h) = - \int _{\partial \Omega} G^h E_{\dashv}^h \ \mathrm{d}\Gamma\quad \quad \forall G^h \;.
\end{equation}
where $F^h$ is expanded as in (\ref{eq:F_expansion}).
Using (\ref{eq:curlF_expansion}) for $\mbox{\bf{curl}}\ F^h$, and (\ref{eq:M0}) and (\ref{eq:M1}) for the mass matrices, we can write the left hand side  of (\ref{eq:var_neum}) as,
\begin{eqnarray}
(\mbox{\bf{curl}}\ G^h, \mbox{\bf{curl}}\ F^h) & + & (G^h, F^h) \nonumber \\
 &=& \mathsf{G}^{\top}
{\mathbb{E}^{1,0}}^{\top}
\mathbb{M}^{(1)} \mathbb{E}^{1,0} \mathsf{F} + \mathsf{G}^{\top} \mathbb{M}^{(0)} \mathsf{F} \;. \label{eq:neu_lhs}
\end{eqnarray}
The boundary conditions are prescribed on the right hand side with the help of duality pairing
\begin{equation} \label{eq:neu_rhs}
\int _{\partial \Omega} G^h E_{\dashv}\ \mathrm{d} \Gamma = G^{\top} \mathbb{T}^{\top} \hat{\mathsf{E}}_{\dashv}^h \;.
\end{equation}
Combining, (\ref{eq:neu_lhs}) and (\ref{eq:neu_rhs}) in (\ref{eq:var_neum}), we have,
\begin{equation}
{\mathbb{E}^{1,0}}^{\top} \mathbb{M}^{(1)} \mathbb{E}^{1,0} \mathsf{F} + \mathbb{M}^{(0)} \mathsf{F} = -\mathbb{T}^{\top} \hat{\mathsf{E}}_{\dashv} \;.
\label{eq:discrete_crulF}
\end{equation}

\subsection{The Dirichlet problem}
The variational formulation for the Dirichlet problem (\ref{eq:curlE}) is given by: Find $\widetilde{\bm{E}}^h$ such that,
\begin{equation}
(\mbox{curl}\ \widetilde{\bm{G}}^h, \mbox{curl}\ \widetilde{\bm{E}}^h) + (\widetilde{\bm{G}}^h,\widetilde{\bm{E}}^h) = -\int _{\partial \Omega} \mbox{curl} \widetilde{\bm{G}}^h \hat{E}_{\dashv}\ \mathrm{d}\Gamma  \quad \quad \forall \widetilde{\bm{G}}^h \;.
\label{eq:weakE}
\end{equation}
Here, $\widetilde{\bm{E}}^h$ is expanded in terms of dual polynomials
\begin{equation}
\widetilde{\bm{E}}^h (\xi, \eta) =  \widetilde{\Psi}^{(1)}(\xi, \eta)  \widetilde{\mathsf{E}} \;.
\end{equation}
Then the weak formulation (\ref{eq:weakE}) can be written as
\begin{eqnarray}
(\mbox{curl}\ \widetilde{\bm{G}}^h, \mbox{curl}\ \widetilde{\bm{E}}^h) & + & (\widetilde{\bm{G}}^h,\widetilde{\bm{E}}^h) \nonumber \\
  & = & \widetilde{\mathsf{G}}^T {\mathbb{E}^{1,0}} \widetilde{\mathbb{M}}^{(2)} {\mathbb{E}^{1,0}}^T \widetilde{\mathsf{E}} + \widetilde{\mathsf{G}} ^T\widetilde{\mathbb{M}}^{(1)} \widetilde{\mathsf{E}} \;,
  \label{eq:GE_Hcurl_inner}
\end{eqnarray}
and
\begin{equation}
\int _{\partial \Omega} \mbox{curl} \widetilde{\bm{G}}^h \hat{E}_{\dashv}\ \mathrm{d}\Gamma  = \widetilde{\mathsf{G}}^T \mathbb{E}^{1,0} \widetilde{\mathbb{M}}^{(2)} \mathbb{T}^T \hat{\mathsf{E}}_{\dashv} \;.
\end{equation}
So the weak formulation (\ref{eq:weakE}) can be written as
\begin{equation}
{\mathbb{E}^{1,0}} \widetilde{\mathbb{M}}^{(2)} {\mathbb{E}^{1,0}}^T \widetilde{\mathsf{E}} + \widetilde{\mathbb{M}}^{(1)} \widetilde{\mathsf{E}} = - \mathbb{E}^{1,0} \widetilde{\mathbb{M}}^{(2)} \mathbb{T}^T \hat{\mathsf{E}}_{\dashv}  \;.
\label{eq:discrete_crulE}
\end{equation}
\subsection{The equivalence condition}\label{sec:equivalence}
In this part, we prove that the two discrete formulations (\ref{eq:discrete_crulF}) and (\ref{eq:discrete_crulE}) are related by the discrete relation $\bm{E}^h = \mbox{\bf{curl}}F^h$, which, in terms of the degrees of freedom is equivalent to $\widetilde{\mathsf{E}}=\mathbb{M}^{(1)}\mathbb{E}^{1,0}\mathsf{F}$.
If we substitute this in the left hand side of (\ref{eq:discrete_crulE}) we get
\begin{equation}
{\mathbb{E}^{1,0}} \widetilde{\mathbb{M}}^{(2)} {\mathbb{E}^{1,0}}^T \widetilde{\mathrm{E}} + \widetilde{\mathbb{M}}^{(1)} \widetilde{\mathrm{E}} = {\mathbb{E}^{1,0}} \widetilde{\mathbb{M}}^{(2)} {\mathbb{E}^{1,0}}^T \mathbb{M}^{(1)} \mathbb{E}^{1,0} \mathrm{F} + \widetilde{\mathbb{M}}^{(1)} \mathbb{M}^{(1)} \mathbb{E}^{1,0} \mathrm{F} \;.
\end{equation}
Then we use (\ref{eq:discrete_crulF}) in the first term on the right hand side and use the fact that $\widetilde{\mathbb{M}}^{(1)}\mathbb{M}^{(1)} = \mathbb{I}$ for the second term on the right hand side to get,
\begin{equation}
{\mathbb{E}^{1,0}} \widetilde{\mathbb{M}}^{(2)} {\mathbb{E}^{1,0}}^T \widetilde{\mathrm{E}} + \widetilde{\mathbb{M}}^{(1)} \widetilde{\mathrm{E}} = {\mathbb{E}^{1,0}} \widetilde{\mathbb{M}}^{(2)} (-\mathbb{T}^T \hat{\mathrm{E}}_\dashv - \mathbb{M}^{(0)} \mathrm{F}) + \mathbb{E}^{1,0} \mathrm{F} \;.
\end{equation}
A further simplification of the bracket terms and using, $\widetilde{\mathbb{M}}^{(2)}\mathbb{M}^{(0)} = \mathbb{I}$, we get,
\begin{equation}
{\mathbb{E}^{1,0}} \widetilde{\mathbb{M}}^{(2)} {\mathbb{E}^{1,0}}^T \widetilde{\mathrm{E}} + \widetilde{\mathbb{M}}^{(1)} \widetilde{\mathrm{E}} = -\mathbb{E}^{1,0} \widetilde{\mathbb{M}}^{(2)} \mathbb{T}^T \hat{\mathrm{E}}_\dashv \;,
\end{equation}
which shows that $\widetilde{\mathsf{E}}$ satisfies (\ref{eq:discrete_crulE}), as required.
\subsection{Equality of norms}\label{sec:equivalent_norms}
The degrees of freedom of the curl of $\bm{E}^h$ are given by (\ref{eq:dirivative_dual_E}) and the associated basis functions in which these degrees are expanded are given by $\widetilde{\Psi}^{(2)}$, (\ref{eq:dual_volume_basis}).
Therefore we have
\begin{eqnarray*}
\| \bm{E}^h\|_{H(\mbox{curl})}^2 & = & \left (\widetilde{\mathsf{E}}^T \mathbb{E}^{1,0} +  \hat{\mathsf{E}}_{\dashv}^T \mathbb{T} \right ) \widetilde{\mathbb{M}}^{(2)} \left (\mathbb{T}^T {\hat{\mathsf{E}}_{\dashv}} + {\mathbb{E}^{1,0}}^T \widetilde{\mathsf{E}} \right ) + \widetilde{\mathsf{E}} ^T\widetilde{\mathbb{M}}^{(1)} \widetilde{\mathsf{E}} \\
 & \stackrel{\widetilde{\mathsf{E}}=\mathbb{M}^{(1)}\mathbb{E}^{1,0}\mathsf{F}}{=} & \left (\mathsf{F}^T {\mathbb{E}^{1,0}}^T \mathbb{M}^{(1)}\mathbb{E}^{1,0} +  \hat{\mathsf{E}}_{\dashv}^T \mathbb{T} \right ) \widetilde{\mathbb{M}}^{(2)} \left ( \mathbb{T}^T {\hat{\mathsf{E}}_{\dashv}} + {\mathbb{E}^{1,0}}^T \mathbb{M}^{(1)} \mathbb{E}^{1,0} \mathsf{F} \right ) + \mathsf{F}^T {\mathbb{E}^{1,0}}^T \mathbb{M}^{(1)} \widetilde{\mathbb{M}}^{(1)} \mathbb{M}^{(1)} \mathbb{E}^{1,0} \mathsf{F} \\
 & \stackrel{(\ref{eq:discrete_crulF})}{=} &  \mathsf{F}^T \mathbb{M}^{(0)} \widetilde{\mathbb{M}}^{(2)} \mathbb{M}^{(0)} \mathsf{F} + \mathsf{F}^T {\mathbb{E}^{1,0}}^T \mathbb{M}^{(1)} \mathbb{E}^{1,0} \mathsf{F} \\
 & \stackrel{\widetilde{\mathbb{M}}^{(2)} \mathbb{M}^{(0)}=\mathbb{I}}{=} &  \mathsf{F}^T \mathbb{M}^{(0)} \mathsf{F} + \mathsf{F}^T {\mathbb{E}^{1,0}}^T \mathbb{M}^{(1)} \mathbb{E}^{1,0} \mathsf{F} \\
 & = &   \| F^h\|_{H(\mbox{\bf{curl}})}^2 \;,
\end{eqnarray*}
which demonstrates that also in the finite dimensional setting the norms of $\bm{E}^h$ and $F^h$ are the same.

\section{Test case} \label{sec:Results}
In this section we show the results for spectral element approximations of (\ref{eq:curlE}) and (\ref{eq:curlF}) for $\Omega \in [-1,1]^2$, using one spectral element.

We choose an exact solution for scalar $F$, and the corresponding vector field $\bm{E}$, given by,
\begin{equation}
F_{ex} = e^x + e^y \;; \quad \quad \bm{E}_{ex} = (e^y, - e^x) \;.
\label{eq:F_and_E_exact}
\end{equation}
The problem (\ref{eq:curlF}) is discretized using a primal representation where we prescribe the Neumann boundary condition,
\begin{equation}
\hat{E}_{\dashv} = \bm{n} \times (\mbox{\bf{curl}}\ F_{ex})  \;.
\end{equation}
For the problem (\ref{eq:curlE}) we use a dual representation where we prescribe the Dirichlet boundary conditions,
\begin{equation}
\hat{E}_{\dashv} = \bm{n} \times \bm{E}_{ex}  \;.
\end{equation}

\begin{figure}
{\includegraphics[width = 0.45\textwidth]{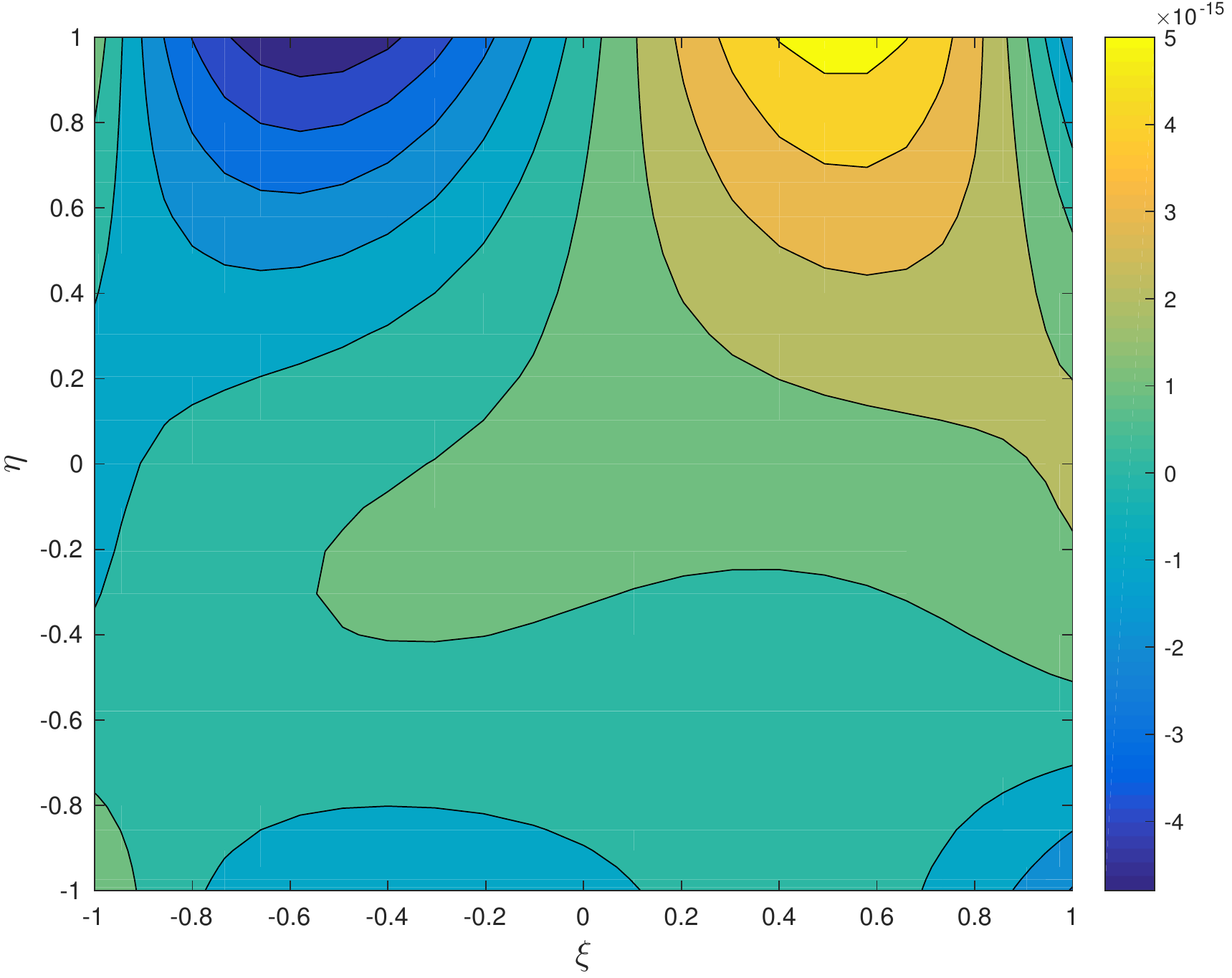}} \quad \quad
{\includegraphics[width = 0.45\textwidth]{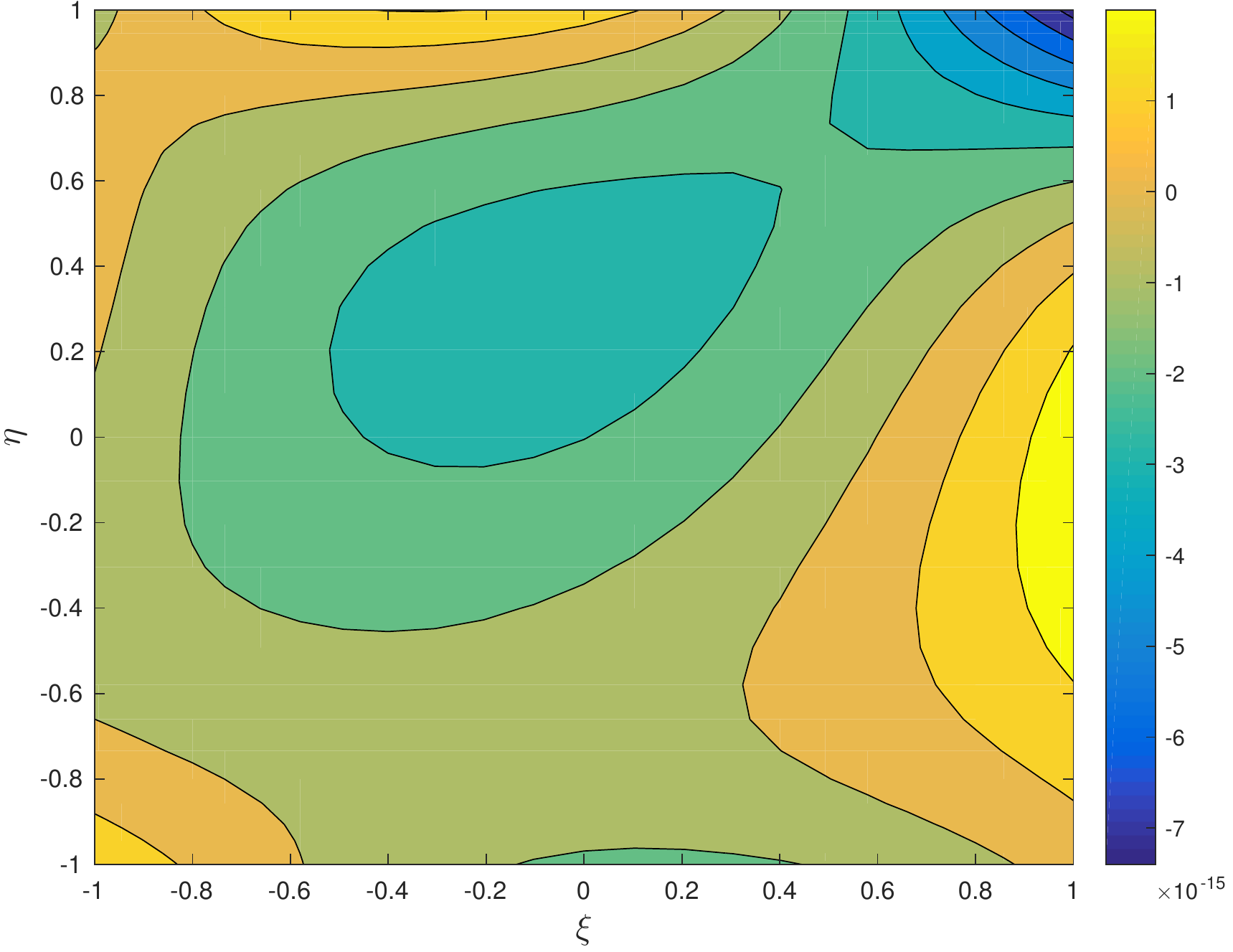}}
\caption{Difference in vector field $\bm{E}^h$ and vector field $(\mbox{\textbf{curl}}\ F^h)$ for $N=3$. \emph{Left}: $\xi$ - component. \emph{Right}: $\eta$ - component.}
\label{fig:equality_c_00}
\end{figure}

In Fig.~\ref{fig:equality_c_00} we show the difference between the $\xi$ and $\eta$ components of vector field $\bm{E}^h$ and $\mbox{\bf{curl}}\ F^h$, for a low order spectral element approximation $N = 3$. The case $N=3$ corresponds to the grid shown in Figure~\ref{fig:F_grid}.
Here, we choose a very low order approximation to show that the equivalence of duality relation derived in Section~\ref{sec:equivalence} holds true even for low order approxiations.

The difference observed in Fig.~\ref{fig:equality_c_00} is of the order $\mathcal{O}(10^{-16})$; the two discrete vector fields agree up to machine precision.
\begin{table}[htbp]
	\centering
	\caption{Norms $\| \bm{E}^h \|_{H(\mbox{curl})}$ and $\|F^h\|_{H(\mbox{\bf{curl}})}$ as a function of the polynomial degree $N$.}
	\begin{tabular}{lcc}
		 \quad $N$ \quad  & \quad $\left\|F^{h}\right\|_{\mathrm{H(\bf{curl})}}$ \quad & \quad $\left\|\bm{E}^{h}\right\|_{H\mathrm{(curl)}}$ \quad \\
		\hline
1 & 5.62334036 &  5.62334036 \\

2 & 6.28815932 & 6.28815932 \\

3 & 6.32851719 & 6.32851719 \\

4 & 6.32957061 & 6.32957061 \\

5 & 6.32958640 & 6.32958640 \\

6 & 6.32958655 & 6.32958655 \\

7 & 6.32958656 & 6.32958656 \\

8 & 6.32958656 & 6.32958656 \\

9 & 6.32958656 & 6.32958656 \\
		\hline
	\end{tabular}%
	\label{tab:norm}%
\end{table}

From Section~\ref{sec:equivalent_norms} we know that in the continuous setting, the $H(\mbox{curl})$-norm of vector field $\bm{E}$ is equal to the $H(\mbox{\bf{curl}})$-norm of scalar field $F$.
For this test case with exact solution given by (\ref{eq:F_and_E_exact}) we have
\begin{equation}
\parallel F_{ex} \parallel_{H(\mbox{curl})} =
\sqrt{\parallel F_{ex} \parallel_{L^2}^2 + \parallel \mbox{\bf{curl}}\ F_{ex} \parallel_{L^2}^2}
= \sqrt{8(\sinh(2)+\sinh^2(1))} = 6.32958656 \;.
\label{eq:theoretical_value}
\end{equation}
In Table~\ref{tab:norm} we show the calculated value of these discrete norms for increasing order of basis functions.
We observe that the discrete norms are exactly equal to each other and they converge to the theoretical value, (\ref{eq:theoretical_value}).

It is worth emphasizing, that the dual basis functions introduced in Section~\ref{sec:dual_formulation}, enforce strong duality pairing that ensures the equivalence of solution, and thus also the equivalence of norms.
\begin{figure}
	\centering
{\includegraphics[width = 0.65\textwidth]{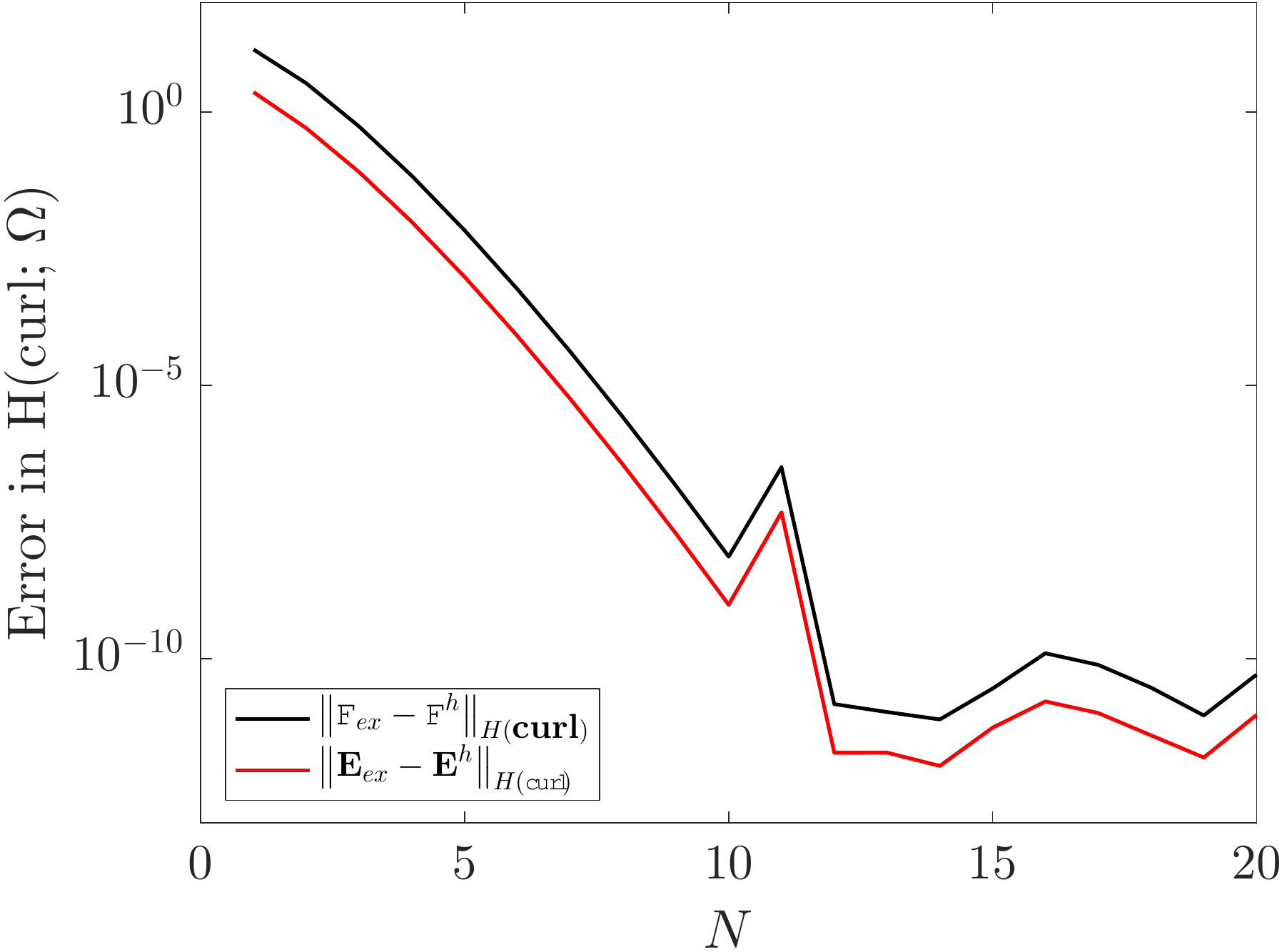}}
\caption{Convergence of $\|F-F^h\|_{H(\mbox{\bf{curl}})}$ and $\|\bm{E}-\bm{E}^h\|_{H(\mbox{curl})}$ as a function of the polynomial order $N$. }
\label{fig:error_convergence_c_00}
\end{figure}

In Fig.~\ref{fig:error_convergence_c_00} we show the convergence of error of $F$ in $H(\mbox{\textbf{curl}};\Omega)$ norm, and the convergence of error of $\bm{E}$ in $\mathrm{H(curl;\Omega)}$ norm.
Both the errors converge exponentially to machine precision level.

\section{Conclusions}
If $\bm{g}$ is the metric tensor in an inner-product space $\mathcal{V}$, the dual degrees of freedom $\alpha_i$ may be related to the expansion coefficients $u^j$ by (\ref{eq:Riesz}), i.e. $\alpha_i = g_{ij} u^j$. Commonly referred to as the {\em raising or lowering of the indices by the metric tensor}. In finite element methods a similar procedure is possible where the role of the metric tensor is played by mass matrices. The two curl-curl problems introduced in Section~\ref{sec:curl_curl_problem} are equivalent in the sense that $F = \mbox{curl}\ \bm{E}$ and the norms of $F$ and $\bm{E}$ are the same. In this paper it is proved that this equivalence continues to hold in finite dimensional spaces, if one of the degrees of freedom is expressed in terms of primal unknowns, $F^h$, and the other in dual degrees of freedom, $\bm{E}^h$. Equivalence of the approximate solutions and their norms is shown in Section~\ref{sec:Discrete_formulations_curl_curl}, while in Section~\ref{sec:Results} this was illustrated for a specific test case.


%
%
%

\end{document}